\documentclass[11pt]{amsart}

\usepackage{xypic}
\xyoption{all}

\usepackage[linktocpage=true]{hyperref}
\newcommand{\myendbibitem}{\relax}

\numberwithin{equation}{section}

\swapnumbers
\newtheorem{thm}[equation]{Theorem}
\newtheorem{prop}[equation]{Proposition}
\newtheorem{lem}[equation]{Lemma}
\newtheorem{cor}[equation]{Corollary}

\theoremstyle{definition}
\newtheorem{defn}[equation]{Definition}
\newtheorem{remark}[equation]{Remark}
\newtheorem{example}[equation]{Example}

\newcommand{\sdp}{\mathbin{{>}\!{\triangleleft}}}

\newcommand{\Aut}{\operatorname{Aut}}

\newcommand{\inv}{^{-1}}

\newcommand{\N}{\operatorname{N}}

\newcommand{\Z}{\operatorname{Z}}

\newcommand{\GL}{{\operatorname{GL}}}

\newcommand{\SL}{{\operatorname{SL}}}

\newcommand{\PGL}{{\operatorname{PGL}}}
\newcommand{\PGLn}{{\operatorname{PGL}_n}}

\newcommand{\Stab}{\operatorname{Stab}}

\newcommand{\Char}{\operatorname{char}} 

\newcommand{\Gal}{\operatorname{Gal}}
\newcommand{\Var}{\operatorname{Var}}
\newcommand{\RMaps}{\operatorname{RMaps}}
\newcommand{\lra}{\longrightarrow}
\newcommand{\Mat}{{\operatorname{M}}}
\newcommand{\M}{\operatorname{M}}
\newcommand{\Mn}{\Mat_n}

\newcommand{\Sym}{{\operatorname{Sym}}}
\newcommand{\End}{{\operatorname{End}}}

\newcommand{\cA}{\mathcal{A}}

\begin{document}
 
\date{July 30, 2007}

\title[Tame actions --- July 30, 2007]
   {Tame group actions on central simple algebras \\ \  \\
 July 30, 2007}

\author{Z. Reichstein}
\address{Department of Mathematics, University of British Columbia,
       Vancouver, BC V6T 1Z2, Canada}
\email{\href{mailto:reichste@math.ubc.ca}%
            {reichste@math.ubc.ca}}
\urladdr{\href{http://www.math.ubc.ca/~reichst}%
              {www.math.ubc.ca/$\stackrel{\sim}{\phantom{.}}$reichst}}

\thanks{Z. Reichstein was supported in part by an NSERC research grant}

\author{N. Vonessen}
\address{Department of Mathematical sciences, University of Montana,
  Missoula, MT 59812-0864, USA}
\email{\href{mailto:Nikolaus.Vonessen@umontana.edu}%
            {Nikolaus.Vonessen@umontana.edu}}
\urladdr{\href{http://www.math.umt.edu/vonessen}{www.math.umt.edu/vonessen}}




\begin{abstract}
  We study actions of linear algebraic groups on finite-di\-men\-sional
  central simple algebras. We describe the fixed algebra for a broad
  class of such actions.
\end{abstract}

\maketitle

\section{Introduction}

Let $k$ be an algebraically closed base field of characteristic zero,
$K/k$ be a finitely generated field extension, $A/K$ be a
finite-dimensional central simple algebra and $G$ be a linear
algebraic group defined over $k$.  Suppose $G$ (more precisely, the
group of $k$-points of $G$) acts on $A$ by $k$-algebra
automorphisms. We would like to describe the fixed ring $A^G$ for this
action, both in terms of its intrinsic structure 
(for example, is it simple?)
and in terms of the way it is embedded in $A$ 
(for example, is the center of $A^G$ contained in $K$?).

The goal of this paper is to give detailed answers to
these and related questions in terms of geometric data, 
for a particular type of $G$-action on $A$ which we call {\em tame}. 
We give a precise definition of tame
actions in Section~\ref{sect.definitions}, and we show in
Section~\ref{sect.sufficient} that many naturally occurring actions
belong to this class.  To every tame action one associates a subgroup
$S$ of $\PGLn$ (defined up to conjugacy), which we call the {\em
  associated stabilizer}.  If the associated stabilizer $S$ is
reductive, we say that the action is {\em very tame}.

To state our first main result, let $K$ be a finitely generated field
extension of $k$.  As usual, we shall say that an algebra $A/K$ is a
{\em form} of $A_0/k$ if the two become isomorphic after extension of
scalars, i.e., $A \otimes_{K} L \simeq A_0 \otimes_k L$ for some field
extension $k \subset K \subset L$.
 
\begin{thm} \label{thm2}
  Consider a tame action of a linear algebraic group $G$ on a central
  simple algebra $A$ of degree $n$ with center~$K$, with associated
  stabilizer subgroup $S \subset \PGLn$.  Then the fixed algebra
  $A^G/K^G$ is a form of $\Mn^S/k$.
\end{thm}

Here $S \subset \PGLn$ is acting on the matrix algebra $\Mn$ by
conjugation, and $\Mn^S=(\Mn)^S$ denotes the fixed subalgebra for this
action.

In the case where the $G$-action on $A$ is very tame, one can
describe $\Mn^S$ (and thus $A^G$) more explicitly.
In this case the associated stabilizer subgroup $S \subset \PGLn$ 
is, by definition, reductive, and hence so is 
the preimage $S^*$ of $S$ in $\GL_n$. 
We shall view the inclusion $\phi \colon S^* \hookrightarrow \GL_n$ 
as an $n$-dimensional linear representation of $S^*$.

 \begin{cor} \label{cor.semisimple}
   Consider a very tame action of a linear algebraic group $G$ on a
   central simple algebra $A$ of degree $n$ with center~$K$, and let
   $\phi$ be as above. Then
 \begin{enumerate}
 \item[(a)] $A^G$ is semisimple.
 \item[(b)] $A^G$ is commutative if and only if $\phi$ is multiplicity-free
 {\upshape(}i.e., a sum of distinct irreducible
 representations{\upshape)}. 
 \item[(c)]$A^G \subset K$ if and only if $\phi$ is irreducible.
 \item[(d)]$\Z(A^G) \subset K$ if and only if 
  $\phi$ is a power of an irreducible representation.
 \end{enumerate}
 \end{cor}
For a more detailed description of $A^G$, in terms of 
the irreducible decomposition of $\phi$, see Section~\ref{sect.very-tame}.

We now return to the setting of Theorem~\ref{thm2}. Let $G$ be 
a linear algebraic group, $A/K$ be a central simple algebra of degree $n$ and
$\psi$ be a tame action of $G$ on $A$, with associated stabilizer 
$S \subset \PGLn$.
Since Theorem~\ref{thm2} asserts that $A^G/K^G$ is a form of $\Mn^S/k$,
it is natural to ask which forms of $\Mn^S/k$ can occur in this way.
Note that here $n$ and $S \subset \PGLn$ are fixed,
but we are allowing $A$, $G$ and $\psi$ to vary.

Our answer to this question will be stated in terms
of Galois cohomology. 
If $G$ is an algebraic group defined over $k$, and
$K/k$ is a field extension, we write $H^1(K, G)$
for the Galois cohomology set 
$H^1(\Gal(K), G(K^{sep}))$, where $K^{sep}$ is 
the separable closure of $K$ and $\Gal(K) = \Gal(K^{sep}/K)$
acts on $G(K^{sep})$ in the natural way; cf.~\cite[Chapter II]{serre-gc}.

For any field extension
$F/k$, the isomorphism classes of the $F$-forms of a $k$-algebra $R$ are in 
a natural 1-1 correspondence with 
$H^1(F, \Aut(R))$; see~\cite[\S 29.B]{boi} or
\cite[Proposition X.4]{serre-lf}.  Here, $\Aut(R)$ is the
algebraic group of algebra automorphisms of $R$. (Recall that 
we are assuming throughout that $\Char(k) = 0$.)  If $B$ is an $F$-form of 
$R$, we will write $[B]$ for the element of $H^1(F, \Aut(R))$ 
corresponding to $B$. 

The normalizer $N = N_{\PGLn}(S)$ naturally acts on the $k$-algebra
$\Mn^S$ (by conjugation). Since $S$ acts trivially, this action 
gives rise to a homomorphism $\alpha \colon N/S \to \Aut(\Mn^S)$ 
of algebraic groups.

\begin{thm} \label{thm.gal-cohomology}
Let $F/k$ be a finitely generated field extension, $S \subset \PGLn$ 
be a closed subgroup and $B$ be an $F$-form of $\Mn^S/k$. 
Then the following are equivalent.
\begin{enumerate}
\item[(a)]There exists a tame action of a linear algebraic group $G$
  on a central simple algebra $A$ of degree $n$ with center $K$, with
  associated stabilizer $S$, such that $A^G/K^G$ is isomorphic to
  $B/F$ {\upshape(}over $k${\upshape)}.
\item[(b)]$[B]$ lies in the image of the natural map 
  \[\alpha_* \colon H^1(F, N/S) \lra H^1(F, \Aut(\Mn^S))\,.\]
\end{enumerate}
\end{thm}

The rest of this paper is structured as follows. In
Section~\ref{sect.definitions}, we introduce tame and very tame
actions.  In Section~\ref{sect.prel.tame} we prove a preparatory lemma
and explain how it is used later on. Sections~\ref{sect.sufficient}
and~\ref{sect.wild} are intended to motivate the notions of tame and
very tame actions. In Sections~\ref{sect.thm2}, \ref{sect.very-tame}
and~\ref{sect.thm.gal-cohomology} we prove Theorem~\ref{thm2},
Corollary~\ref{cor.semisimple} and Theorem~\ref{thm.gal-cohomology}
along with some generalizations.  In the last section we illustrate
Theorem~\ref{thm.gal-cohomology} with an example. This example shows,
in particular, that every finite-dimensional semisimple algebra (over
a finitely generated field extension of $k$) occurs as $A^G/K^G$ for
some very tame action of a linear algebraic group $G$ on some
central simple algebra $A$ with center $K$; see Corollary~\ref{cor3}.

\subsection*{Notation and Terminology}
We shall work over an algebraically closed base field $k$ of
characteristic zero. All fields will be assumed to be finitely
generated extensions of $k$. By a variety we shall mean a reduced 
$k$-scheme of finite type. All varieties, algebraic groups, 
group actions, morphisms, rational maps, etc., are assumed to 
be defined over $k$.  By a point of an algebraic variety we
shall always mean a $k$-point. We will often identify 
an algebraic group $G$ with its group $G(k)$ of $k$-points. 
Every central simple algebra is assumed (by definition) to 
be finite-dimensional over its center.  
 
We denote by $\Mn$ the algebra of $n\times n$ matrices over $k$.  If
$S$ is a subgroup of $\PGLn$, we denote by $\Mn^S$ the fixed algebra
for the conjugation action of $S$ on $\Mn$.  Throughout, $G$ will be a
linear algebraic group.  We shall refer to an algebraic variety $X$
endowed with a regular $G$-action as a $G$-variety.  We will say that
$X$ (or the $G$-action on $X$) is generically free if $\Stab_G(x) = \{
1 \}$ for $x \in X$ in general position.

\subsection*{Acknowledgment} We are grateful to the referee 
for a careful reading of the paper and constructive comments.

\section{Tame actions on central simple algebras}
\label{sect.definitions}

Tame actions are a special class of {\it geometric} actions.  To
review the latter concept, we need to recall the relationship
between central simple algebras of degree~$n$ and irreducible,
generically free $\PGLn$-varieties.

If $X$ is a $\PGLn$-variety, we shall denote the algebra of
$\PGLn$-invariant rational maps $X \dasharrow \Mn$ by
$\RMaps_{\PGLn}(X, \Mn)$. The algebra structure in $\RMaps_{\PGLn}(X,
\Mn)$ is naturally induced from $\Mn$, i.e., given $a, b \colon X
\dasharrow \Mn$, we define $a + b$ and $ab$ by $(a+b)(x) = a(x) +
b(x)$ and $ab(x) = a(x) b(x)$ for $x \in X$ in general position. If
the $\PGLn$-action on $X$ is generically free, then $\RMaps_{\PGLn}(X,
\Mn)$ is a central simple algebra of degree $n$, with center
$k(X)^{\PGLn}$; cf.\ \cite[Lemmas 8.5 and 9.1]{reichstein}.  In this
case we will sometimes denote this algebra by $k_n(X)$. As this
notation suggests, $k_n(X)$ shares many properties with the function
field $k(X)$ (see~\cite{rv2}); in particular, the two coincide if $n =
1$.  Under our hypotheses, every central simple algebra of degree $n$
is isomorphic to $k_n(X)$ for some irreducible generically free
$\PGLn$-variety $X$; see~\cite[Theorem 1.2]{rv2}.

Suppose that $X$ is a $G \times \PGLn$-variety, and that the
$\PGLn$-action on $X$ is generically free. Then the $G$-action on $X$
naturally induces a $G$-action on $k_n(X)$; see~\cite{rv3}.
Following~\cite{rv3} we will say that the action of an algebraic group
$G$ on a central simple algebra $A/K$ is {\em geometric} if $A$
is $G$-equivariantly isomorphic to $k_n(X)$ for some $G \times
\PGLn$-variety $X$ as above.
The $G \times \PGLn$-variety $X$ is then
called the {\em associated variety} for the $G$-action on $A$. This
associated variety is unique up to birational isomorphism of $G \times
\PGLn$-varieties, see \cite[Corollary 3.2]{rv3}. 

Note that the center of $k_n(X)$ is $k(X/\PGLn)$.  
So a central simple algebra~$A/K$ cannot admit a geometric action 
unless $K$ is finitely generated over $k$.  For this reason 
we will only consider central simple algebras whose centers 
are finitely generated extensions of $k$ throughout this paper.

The class of geometric actions on central simple algebras 
is rather broad; it includes, in particular, all 
``algebraic'' actions; see \cite[Theorem 5.3]{rv3}.
Roughly speaking, an action of $G$ on $A$ is algebraic 
if there is a $G$-invariant order $R$ in $A$ such that $G$ 
acts regularly (``rationally'') on $R$. For details, see \cite{rv3}.

Consider a geometric action $\psi$ of a linear algebraic group $G$ on
a central simple algebra $A$, with associated $G \times \PGLn$-variety
$X$.  Given $x \in X$, denote by $S_x$ the projection of $\Stab_{G
  \times \PGLn}(x) \subset G \times \PGLn$ to $\PGLn$.

\begin{defn} \label{def.tame}
  (a) We shall call the action $\psi$ {\em tame} if there is a dense
  open subset $U \subset X$ such that $S_x$ is conjugate to $S_y$ in
  $\PGLn$ for every $x, y \in U$. In this case we shall refer to $S_x$
  ($x \in U$) as the {\em associated stabilizer}%
  \footnote{For $x\in X$ in general position, $S_x$ is the stabilizer 
  of the image of $x$ in the rational quotient $X/G$, 
 see Remark~\ref{rem.prel1}(b).}
  for $\psi$; the associated stabilizer is defined up to conjugacy.
  (b) We shall call the action $\psi$ {\em very tame} if it is tame
  and if its associated stabilizer is reductive.
\end{defn}

See Lemma~\ref{lem.tame} for an alternative description of tame and very
tame actions.  There exist non-tame geometric actions, and tame
actions that are not very tame (see Proposition~\ref{prop.wild}).
However, informally speaking, many interesting actions are geometric
and very tame; in particular, all geometric actions on division 
algebras are very tame. For details, see Section~\ref{sect.sufficient}.

\section{An exact sequence of stabilizer groups}
\label{sect.prel.tame}

\begin{lem} \label{prel1.tame} 
  Let $\Gamma$ be a linear algebraic group, $N$ a closed normal subgroup,
  and $\pi \colon \Gamma \lra \Gamma/N$  the natural projection.
  Suppose $X$ is a $\Gamma$-variety, $Y$ is a birational model for
  $X/N$ {\upshape(}as a $\Gamma/N$-variety{\upshape)}, and $q \colon X
  \dasharrow Y$ 
  is the rational quotient map. Then
  for $x \in X$ in general position, the sequence of
  stabilizer subgroups
  \[ \qquad 1 \lra \Stab_N(x) \lra \Stab_{\Gamma}(x) \stackrel{\pi}{\lra} 
  \Stab_{\Gamma/N}(q(x)) \lra 1 \]
  is exact.
\end{lem}

\begin{proof}
  The sequence is clearly exact at $\Stab_N(x)$ and
  $\Stab_{\Gamma}(x)$, so we only need to show that the map \[
  \Stab_{\Gamma}(x) {\lra} \Stab_{\Gamma/N}(q(x)) \, , \]
which we denote by $\pi$, is
  surjective.  By a theorem of Rosenlicht 
(see \cite{rosenlicht1}, \cite{rosenlicht2} or \cite[Proposition 2.5]{pv}) 
  there is
  a dense $\Gamma$-invariant open subset $X_0 \subset X$ such that
  $q^{-1}(q(x)) = N \cdot x$ for every $x \in X_0$.  Thus for $x \in
  X_0$, $g \in \Gamma$,
  \begin{align*}
  \text{$\pi(g) \in G/N$ stabilizes $q(x)$}
    &\Longleftrightarrow
           \text{$g(x) \in N \cdot x$}\\
    &\Longleftrightarrow
           \text{$n^{-1}g\in\Stab_{\Gamma}(x)$ for some $n\in N$} \\
    &\Longleftrightarrow
           \text{$\pi(g) = \pi(n^{-1} g)\in\pi(\Stab_{\Gamma}(x))$}\,,
  \end{align*}
  as claimed.  
\end{proof}

\begin{remark} \label{rem.prel1} In the sequel we repeatedly
apply Lemma~\ref{prel1.tame} in a situation, where
$\Gamma = G \times \PGLn$ and $X$ is a $\Gamma$-variety 
with a generically free $\PGLn$-action. There are two natural choices 
of $N$ here, $N = G$ and $N = \PGLn$ (or, more precisely,
$N = G \times \{ 1 \}$ and $N = \{ 1 \} \times \PGLn$). Thus we will 
consider two rational quotient maps:
\[  \xymatrix{    &  X \ar@{-->}[ld] \ar@{-->}[rd] & \cr 
       Y = X/G  &                                & Z = X/\PGLn \, . }
\]

(a) Taking $N = \PGLn$ in Lemma~\ref{prel1.tame} and
remembering that the $\PGLn$-action on $X$ is generically free,
we see that 
the natural projection $G \times \PGLn \lra G$ restricts to an isomorphism 
\[ \Stab_{G \times \PGLn}(x) \stackrel{\simeq}{\lra}
\Stab_{G}(z) \, , \] 
where $x \in X$ is in general position and $z$ is the image 
of $x$ in $X/\PGLn$; cf.~\cite[Lemma 7.1]{rv3}.

\smallskip
(b) On the other hand, taking $N = G$ and
denoting by $\pi$ the natural projection $G \times \PGLn \lra \PGLn$
we obtain 
\[ \Stab_{\PGLn}(y) = \pi(\Stab_{G \times \PGLn}(x)) \simeq 
\Stab_{G \times \PGLn}(x) / \Stab_{G}(x) \, , \] 
where $x \in X$ is in general position and $y$ is the image 
of $x$ in $X/G$, cf.\ \cite[Lemma~2.3]{rv4}. 
\end{remark}

\section{An alternative description of tame and very tame actions, and
  some examples}
\label{sect.sufficient}

It follows from Remark~\ref{rem.prel1}(b) that
if $\psi$ is a geometric $G$-action on a central simple
algebra $A$, with associated $G \times \PGLn$-variety $X$, then $\psi$
is tame if and only if the $\PGLn$-action on $X/G$ has a {\em
stabilizer in general position} $S \subset \PGLn$,
in the sense of \cite[p. 228]{pv}. Recall that this means that
$\Stab_{\PGL_n}(y)$ is conjugate to $S$ for a general point $y \in X/G$,
and that $S$ is only well defined up to conjugacy in $\PGLn$.
Recall also that if $\psi$ is tame, we called $S$ the {\em associated 
stabilizer} for $\psi$; see Definition~\ref{def.tame}.  

Note that the stabilizer in general position is known to exist under
certain mild assumptions on the action; cf.~\cite[\S7]{pv}. Informally 
speaking, this means that ``most'' geometric actions on 
central simple algebras are tame; we illustrate this point by 
the examples below.  In particular, if $\Stab_\PGLn(y)$ 
is reductive for general $y \in X/G$ then, by a theorem 
of Richardson~\cite[Theorem 9.3.1]{richardson}, $X/G$ has a stabilizer 
in general position and thus the action $\psi$ is very tame.
For later use, we record some of these remarks as a lemma.

\begin{lem} \label{lem.tame}
  Consider a geometric action of a linear algebraic group $G$ on a
  central simple algebra~$A$ with associated
  $G\times\PGLn$-variety~$X$.
  \begin{enumerate}
  \item[(a)]The action of $G$ on $A$ is tame if and only if the
    $\PGLn$-action on $X/G$ has a stabilizer in general position.
  \item[(b)]The action of $G$ on $A$ is very tame if and only if
    $\Stab_\PGLn(y)$ is reductive for $y\in X/G$ in general position.
    \qed
  \end{enumerate}
\end{lem}

\begin{example} \label{ex.tame1} 
A geometric action $\psi$ on a central simple algebra $A$ is very tame
if

\begin{itemize}
 \item[(a)]$A$ is a division algebra, or
 \item[(b)] the connected component $G^0$ of $G$ is a torus (this
 includes the case where $G$ is finite), or
\item[(c)] $\Stab_G(z)$ is reductive for $z\in X/\PGLn$ in general position, or
\item[(d)] $G$ is reductive and the $G$-action on $A$ is inner, or
 \item[(e)] $A$ has a $G$-invariant maximal \'etale subalgebra.
\end{itemize}
\end{example}

Here by an \'etale subalgebra of $A$ we mean a subalgebra of $A$ which
is a (finite) direct sum of finite field extensions of the center
of~$A$ (since we are working in characteristic zero, these field
extensions are necessarily separable).

\begin{proof}
  Using Richardson's theorem (see the remarks before
  Lemma~\ref{lem.tame}), part~(a) follows from \cite[Section 3]{rv4},
  so we will focus on the other parts.  Let $X$ be the associated $G
  \times \PGLn$-variety, let $x \in X$ be a point in general position,
  and let $y$ and $z$ be the images of $x$ in $Y = X/G$ and $Z =
  X/\PGLn$, respectively, as in Remark~\ref{rem.prel1}. By
  Lemma~\ref{lem.tame}, it suffices to show that $\Stab_{\PGLn}(y)$ is
  reductive for $y \in Y$ in general position.

\smallskip
  (b) For $x \in X$ in general position, $\Stab_{G \times \PGLn}(x)$ 
is isomorphic to a closed subgroup of $G$ and is thus reductive;
   see Remark~\ref{rem.prel1}(a). Hence,
\[ \Stab_{\PGLn}(y) \simeq 
\Stab_{G \times \PGLn}(x) / \Stab_{G}(x) \] 
is also reductive for $y \in Y$ in general position; 
cf.~Remark~\ref{rem.prel1}(b).

\smallskip
(c) Let $x\in X$ in general position.  
  Since $\Stab_{G \times \PGLn}(x) \simeq
  \Stab_G(z)$, it is a reductive group, and so is its homomorphic
  image $\Stab_{\PGLn}(y)$.  The desired conclusion follows from
  Lemma~\ref{lem.tame}.

  \smallskip
(d) Here $G$ acts trivially on $\Z(A)=k(X/\PGL_n)$ and hence, on
$X/\PGL_n$.  Thus $\Stab_G(z) = G$ is reductive for every $z \in X/\PGL_n$,
and part (c) applies.

  \smallskip
  (e) By \cite[Theorem 1.5(b)]{rv3} for $x\in X$ in
  general position there exists a
  maximal torus $T \subset \PGLn$ (depending on $x$) such that
  $\Stab_{G \times \PGLn}(x) \subset G \times \N(T)$,
  where $\N(T)$ denotes the normalizer of $T$ in $\PGLn$.  Thus
  the image $\pi(\Stab_{G \times \PGLn}(x))$ of 
  $\Stab_{G \times \PGLn}(x)$ under the natural
  projection $\pi \colon G \times \PGLn \lra \PGLn$ consists of
  semisimple elements and, consequently, is reductive. 
  By Remark~\ref{rem.prel1}(b)
  $\Stab_{\PGLn}(y)$ and $\pi(\Stab_{G \times \PGLn}(x))$ 
  are isomorphic. Hence, $\Stab_{\PGLn}(y)$ is reductive, 
  and part (e) follows.
\end{proof}

\section{Examples of wild actions}
\label{sect.wild}

The purpose of this section is to show that Definition~\ref{def.tame} 
is not vacuous. 

\begin{prop} \label{prop.wild}
{\upshape(a)} There exists a geometric action of an algebraic group 
on a central simple algebra which is tame but not very tame.

\smallskip
{\upshape(b)} There exists a geometric action of an algebraic group 
on a central simple algebra which is not tame.
\end{prop}

Our proof of Proposition~\ref{prop.wild} will rely on
the following lemma.

\begin{lem} \label{lem.wild}
Let $G$ be an algebraic group and $H$ a closed subgroup.
Then for every $H$-variety $Y$ {\upshape(}not necessarily generically
free{\upshape)} 
there exists a $G \times H$-variety $X$ such that
\begin{itemize}
\item[(i)]the action of $H=\{1\}\times H$ on $X$ is generically free, and
\item[(ii)]$X/G$ is birationally isomorphic to $Y$ {\upshape(}as an
  $H$-variety{\upshape)}. 
\end{itemize}
\end{lem}

\begin{proof}
Set $X = G \times Y$, with $G$ acting
on the first factor by translations on the left
and $H$ acting by $h  \colon (g, y) \mapsto (g h^{-1}, hy)$. 
The $G$-action and the $H$-action on $X$ commute and thus
define a $G \times H$-action. Conditions (i) and (ii)
are now easy to check.
\end{proof}

\begin{proof}[Proof of Proposition~\ref{prop.wild}] (a)
  In view of Lemma~\ref{lem.tame}, we need to construct a $G \times
  \PGLn$-variety $X$, with $\PGLn$ acting generically freely, and such
  that the $\PGLn$-action on the rational quotient variety $X/G$ has a
  non-reductive stabilizer in general position.  By
  Lemma~\ref{lem.wild} (with $G = H = \PGLn$) it suffices to construct
  a $\PGLn$-variety $Y$ with a non-reductive stabilizer in general
  position. To construct such a $Y$, take any non-reductive subgroup
  $S \subset \PGLn$ and let $Y$ be the homogeneous space $\PGLn/S$.

\smallskip
(b) Arguing as in part (a), we see that 
it suffices to construct a $\PGLn$-variety $Y$ which does not have 
a stabilizer in general position. Richardson~\cite[\S12.4]{richardson} 
showed that (for suitably chosen $n$) there exists 
an $\SL_n$-variety $Z$ with the following property:
for every nonempty Zariski open subset $U \subset Z$ there are infinitely many
stabilizers $\Stab_{\SL_n}(z)$, $z \in U$, with pairwise non-isomorphic
Lie algebras.   Clearly such an $\SL_n$-variety $Z$ cannot have 
a stabilizer in general position. We will now construct 
a $\PGLn$-variety $Y$ with no stabilizer in general position 
by modifying this example. Note that in~\cite{richardson} 
the variety $Z$ is only constructed over the field $\mathbb{C}$ 
of the complex numbers. However, by the Lefschetz principle
this construction can be reproduced over any algebraically 
closed base field $k$ of characteristic zero.

Consider the rational quotient map $\pi \colon Z \dasharrow Y = Z/\mu_n$,
for the action of the center $\mu_n$ of $\SL_n$ on $Z$. Recall that
(after choosing a suitable birational model for $Y$), the
$\SL_n$-action on $Z$ descends to a $\PGLn$-action on $Y$. 
In fact, this is exactly the situation we considered in the setting of 
Lemma~\ref{prel1.tame} (with $\Gamma = \SL_n$ and $N = \mu_n$); 
this lemma tells us that, after replacing $Y$ by a dense 
open subset $Y_0$, and $Z$ by $\pi^{-1}(Y_0)$, 
\[ \Stab_{\PGLn}(y) \simeq \Stab_{\SL_n}(z)/\Stab_{\mu_n}(z) \]
for every $z \in Z$ and $y = \pi(z)$.
Since $\Stab_{\mu_n}(z)$ is a finite group we see that
$\Stab_{\PGLn}(y)$ and $\Stab_{\SL_n}(z)$ have isomorphic Lie algebras.
This shows that for every nonempty open subset $V$ of $Y$ there are infinitely
many stabilizers $\Stab_{\PGLn}(y)$, $y \in V$,
with pairwise non-isomorphic Lie algebras.
In particular, the $\PGLn$-action on $Y$ has 
no stabilizer in general position.
\end{proof}

\section{Homogeneous fiber spaces}
\label{sect.fiber-spaces}

In this section we briefly recall the definition and basic properties
of homogeneous fiber spaces. Let $\varphi\colon H\to G$ be 
a homomorphism of linear algebraic groups (defined over $k$) and
$Z$ be an $H$-variety. The {\it homogeneous fiber space} 
$X = G*_H Z$ is defined as the rational quotient of $G \times
Z$ for the action of $\{1\}\times H$, where $G \times H$ acts on $G
\times Z$ via $(g, h) \cdot (\tilde{g}, z)=(g \tilde{g} \varphi(h)\inv, hz)$.
As usual, we choose a model for this rational quotient, so that the
$G$-action on $G \times Z$ descends to a regular action on this model.
Note that $\{ 1_G \} \times Z$ projects to an $H$-subvariety of $X$
birationally isomorphic to $Z/\ker(\varphi)$.

It is easy to see that $G *_H Z$ is birationally isomorphic
to $G *_{\varphi(H)} (Z/\ker(\varphi))$ as a $G$-variety; thus in carrying
out this construction we may assume that $H$ is a subgroup of $G$ 
and $\varphi$ is the inclusion map. For a more detailed discussion
of homogeneous fiber spaces in this setting and further references,
we refer the reader to \cite[\S3]{rv1}.

We also recall that a $G$-variety $X$ is called {\em primitive} if
$G$ transitively permutes the irreducible components of $X$ or,
equivalently, if $k(X)^G$ is a field; see~\cite[Lemma 2.2(b)]{reichstein}. 

Our goal in this section is to prove
the following variant of~\cite[Lemma 2.17]{lpr}, which
will be repeatedly used in the sequel.

\begin{lem} \label{lem.Y->Z} 
Let $\varphi\colon H\to G$ be 
a homomorphism of linear algebraic groups {\upshape(}defined over
$k${\upshape)}, 
$Z$ be an $H$-variety and $R$ be a finite-dimensional $k$-algebra.  
Let $X$ be the homogeneous fiber space $G *_H Z$.  Assume 
further that $G$ acts
on $R$ via $k$-algebra automorphisms {\upshape(}i.e., via a morphism
$G \lra \Aut(R)$ of algebraic groups{\upshape)}. Then
  \begin{itemize}
  \item[(a)] {\upshape(}cf.~\cite[p. 161]{pv}{\upshape)} $k(X)^G$ and
    $k(Z)^H$ are isomorphic as $k$-algebras. In particular, $X$ is
    primitive if and only if $Z$ is primitive.
  \item[(b)] If $X$ and $Z$ are primitive, the algebras $\RMaps_G(X,
    R)/k(X)^G$ and $\RMaps_H (Z, R)/ k(Z)^H$ are isomorphic.
  \end{itemize}
\end{lem}

Here, as usual, $\RMaps_H(Z, R)$ denotes the $k(Z)^H$-algebra 
of $H$-equivariant rational maps $Z \dasharrow R$.

\begin{proof}
  Let $f \colon X=G *_H Z\dasharrow R$ be a $G$-equivariant rational
  map.  The indeterminacy locus of such a map is a $G$-invariant
  subvariety of $X$. 
  Let $\chi \colon Z \dasharrow X$ be the rational map sending 
  $z \in Z$ to the image 
  of $(1_G,z)\in G\times Z$ in $X = G*_H Z$.
  Since $G \cdot \chi(Z)$ is dense in $X$, $\chi(Z)$ 
  cannot be contained in the indeterminacy locus of $f$. 
  Thus the map 
  \[ \psi \colon \RMaps_G(X, R) \to \RMaps_H (Z, R) \]
  given by $\psi(f)=f\circ\chi$ is a well-defined $k$-algebra homomorphism.
  Since $G \cdot \chi(Z)$ is dense in $X$, $\psi$ is injective.  

  It is easy to see that $\psi$ has an inverse. Indeed, given $H$-equivariant
  rational map $\lambda \colon Z \dasharrow R$, we define a
  $G$-equivariant map $\lambda' \colon G \times Z \dasharrow R$ 
  by $(g, z) \mapsto g \lambda(z)$. Since $\lambda'(g \varphi(h)\inv, h z) =
  \lambda'(g, z)$ for every $g \in G$ and $z \in Z$, and $X=(G \times
  Z)/H$, the universal property of rational quotients tells us that
  $\lambda'$ descends to a $G$-equivariant rational map $f \colon X
  \dasharrow R$. Moreover, by our construction $\psi(f) = f \circ \chi
   = \lambda$, so that the map $\RMaps_H (Z, R) \to \RMaps_G(X, R)$ 
  taking $\lambda$ to $f$ is the inverse of $\psi$.
  
  We are now ready to complete the proof of Lemma~\ref{lem.Y->Z}. 
  To prove part (a), we apply the above construction in the case 
  where $R = k$, with trivial $G$-action. Here $\RMaps_G(X, k) = k(X)^G$,
  $\RMaps_H(Z, k) = k(Z)^H$, and $\psi$ is the desired isomorphism
  of $k$-algebras.  To prove part (b), note that if we identify elements of
  $k(X)^G$ (respectively, $k(Z)^H$) with $G$-equivariant rational maps
  $X \dasharrow R$ (respectively, $H$-equivariant rational maps $Z
  \dasharrow R$) whose image is contained in $k \cdot 1_R$ then
  $\psi$ is an algebra isomorphism between $\RMaps_G(X, R)/ k(X)^G$ 
 and $\RMaps_H (Z, R)/ k(Z)^H$.
\end{proof}

\section{Proof of Theorem~\ref{thm2}}
\label{sect.thm2}

We begin with a simple observation.  Let $X$ be the associated
$G \times \PGLn$-variety for the $G$-action on~$A$.  Let $Y$ be a
$\PGLn$-variety, which is a birational model for $X/G$. Then by the
universal property of rational quotients (see, e.g.,~\cite[\S2.3-2.4]{pv},
\cite[Remark 2.4]{reichstein} or \cite[Remark 6.1]{rv3}),
\begin{equation} \label{e.X->Y}
A^G = \RMaps_{\PGL_n}(X, \Mn)^G \simeq  \RMaps_{\PGLn}(Y, \Mn) \, . 
\end{equation}
Moreover,
\[\Z(A)^G=(k(X)^\PGLn)^G=k(X)^{G\times\PGLn}\simeq k(Y)^\PGLn\,.\]
Thus Theorem~\ref{thm2} is a special case of the following result
(applied to the conjugation action of $\Gamma = \PGLn$ on $R = \Mn$).

\begin{prop} \label{prop.thm2} 
  Let $R$ be a finite-dimensional $k$-algebra, and $\Gamma$ a
  linear algebraic group, acting on $R$ by $k$-algebra
  automorphisms {\upshape(}i.e., via a morphism $\Gamma \lra \Aut(R)$ of
  algebraic groups{\upshape)}.  If a primitive $\Gamma$-variety $Y$ has a
  stabilizer $S \subset \Gamma$ in general position then
  $\RMaps_{\Gamma}(Y, R)/k(Y)^\Gamma$ is a form of $R^S/k$.
\end{prop}

Here we view $\RMaps_{\Gamma}(Y, R)$ as a $K$-algebra, where $K =
k(Y)^{\Gamma}$.  Note that $K$ is a field since $Y$ is a primitive
$\Gamma$-variety (the latter means that $\Gamma$ transitively permutes
the irreducible components of~$Y$). We also recall that 
the stabilizer in general position is only defined up 
to conjugacy in $\Gamma$; we choose a particular subgroup 
$S \subset \Gamma$ in this conjugacy class. Of course, 
replacing $S$ by a conjugate subgroup does not change the
isomorphism type of $R^S/k$. 

\begin{proof} We will proceed in two steps.

\smallskip
{\bf Step 1.} Suppose first that $S = \{ 1 \}$,
i.e., that the $\Gamma$-action on $Y$ is generically free.
Let $Y_0$ be an irreducible component of $Y$, and let $\phi\colon
\Gamma\times Y_0\lra Y$ be the natural $\Gamma$-equivariant map, where
$\Gamma$ acts on itself by left multiplication and trivially on $Y_0$.
Then $\phi$ induces an injective homomorphism
$\RMaps_{\Gamma}(Y, R) \hookrightarrow \RMaps_{\Gamma}(\Gamma \times
Y_0, R)$ via $f \mapsto f \circ \phi$. Clearly
$\Gamma$-equivariant rational maps from $\Gamma\times Y_0$ to $R$
are in 1-1 correspondence with rational maps from $Y_0$ to $R$. In
other words,
\[ \RMaps_{\Gamma}(Y, R) \hookrightarrow 
\RMaps_{\Gamma}(\Gamma\times Y_0, R) \simeq \RMaps(Y_0, R) \simeq R
\otimes_k k(Y_0) \, . \] 
In particular, 
\[ \dim_{k(Y_0)} \RMaps_{\Gamma}(\Gamma \times Y_0, R) = 
\dim_{k(Y_0)}(R \otimes_k k(Y_0)) = \dim_k(R) \, . \]
On the other hand, by~\cite[Lemma 7.4(b)]{reichstein}, 
\[ \dim_K \, \RMaps_{\Gamma}(Y, R) = \dim_k(R) \, . \] 
Now suppose $f_1,\ldots,f_m\in\RMaps_\Gamma(Y,R)$.  Applying
\cite[Lemma 7.4(a)]{reichstein} twice, one sees that the $f_i$ are
linearly independent over $K=k(Y)^\Gamma$ iff the $f_i(y)$ are
$k$-linearly independent for $y\in Y$ in general position iff the
$f_i\circ\phi$ are linearly independent over $k(Y_0)$.  This shows
that the induced map of $k(Y_0)$-algebras
\[ \RMaps_{\Gamma}(Y, R) \otimes_K k(Y_0) \to
\RMaps_\Gamma(\Gamma \times Y_0, R)\]
is injective. Since
$\RMaps_{\Gamma}(Y, R) \otimes_K k(Y_0)$ and $
\RMaps_\Gamma(\Gamma \times Y_0, R)$ have the same dimension 
over $k(Y_0)$, this map is an isomorphism. 
In other words,  $\RMaps_{\Gamma}(Y, R)$
is a $K$-form of $R$. This completes the proof of
Proposition~\ref{prop.thm2} in the case that $S=\{1\}$.

\smallskip 
{\bf Step 2.} 
We will now reduce the general case to the case considered in Step 1.
Let $Z$ be the union of the components of $Y^S$ of maximal dimension,
and let $N = N_{\Gamma}(S)$ be the normalizer of $S$ in $\Gamma$.
Note that $Z$ is an $N$-variety.  
By \cite[Lemma~3.2]{rv1}, 
$Y$ is birationally isomorphic to
$\Gamma*_N Z$, as a $\Gamma$-variety. 
By Lemma~\ref{lem.Y->Z}, $Z$ is primitive since $Y$ is, and
\begin{equation} \label{e.Y->Zrevived}
  \RMaps_{\Gamma}(Y, R) {\simeq} \RMaps_N (Z, R)\
  \text{\ \ and\ \ }
  k(Y)^{\Gamma} \simeq k(Z)^N \, .
\end{equation}
Thus it suffices to show that $\RMaps_N (Z, R)/k(Z)^N$ is a form of
$R^S/k$. 

Recall that $S$ acts trivially on $Z$ (because $Z$ is, by definition,
a subset of $Y^S$); hence, the image of every $N$-equivariant map from
$Z$ to $R$ will actually lie in $R^S$. In other words,
\begin{equation} \label{e.Z}
 \RMaps_N (Z, R) = \RMaps_N(Z, R^S) = \RMaps_{N/S} (Z, R^S) \, . 
\end{equation}
Moreover, $k(Z)^N=k(Z)^{N/S}$.  We now
show that the $N/S$-action on $Z$ is generically free;
the desired conclusion will then follow from the result of Step 1.
Equivalently, we
want to show that $\Stab_N(z) = S$ for $z \in Z$ in general position.
Clearly, $S \subset \Stab_N(z)$ for any $z \in Z \subset Y^S$. To
prove the opposite inclusion, it is enough to show that $S =
\Stab_\Gamma(z)$ for $z \in Z$ in general position.  Indeed, by the
definition of $S$, $\Stab_{\Gamma}(y)$ is conjugate to $S$ for $y \in
Y$ in general position.  Since $\Gamma Z$ is dense in $Y$,
$\Stab_\Gamma(z)$ is conjugate to $S$ for $z \in Z$ in general
position.  Since we know that $S \subset \Stab_\Gamma(z)$ we conclude that
$S = \Stab_\Gamma(z)$ for $z \in Z$ in general position, as desired.
\end{proof}

We conclude this section with a simple corollary of Theorem~\ref{thm2}.

\begin{cor} \label{cor.div} Consider a tame action of $G$ on a central
  simple algebra $A$
with associated stabilizer $S \subset \PGLn$. Let $N$ be the normalizer 
of $S$ in $\PGLn$.  If $A^G$ is a division algebra then $\Mn^{N} = k$.
 \end{cor}
 
\begin{proof} Let $X$ be the associated $G \times \PGLn$-variety.
  Set $Y = X/G$, and denote by $Z$ the union of the irreducible
  components of $Y^S$ of maximal dimension. Recall from the proof of
  Theorem~\ref{thm2} that
  \[ A^G \simeq \RMaps_{\PGLn}(Y, \Mn) \simeq \RMaps_{N}(Z, \Mn^S) \, . \]
  Now observe that if $m \in \Mn^{N}$, then the constant map $f_m
  \colon Z \lra \Mn^S$, given by $f_m(z) = m$ for every $z \in Z$, is
  $N$-equivariant. Thus $f_m$ may be viewed as an element of $A^G$,
  and $m \mapsto f_m$ defines an embedding $\Mn^{N} \hookrightarrow
  A^G$.  Since we are assuming that $A^G$ is a division algebra, this
  implies that $\Mn^{N}$ cannot have zero divisors. On the other hand,
  since $\Mn^{N}$ is a finite-dimensional $k$-algebra and $k$ is
  algebraically closed, this is only possible if $\Mn^{N} = k$.
\end{proof}

\section{Proof of Corollary~\ref{cor.semisimple}}
\label{sect.very-tame}

In this section we spell out what Theorem~\ref{thm2} says in the
case where the $G$-action on $A$ is very tame. 
Here the associated stabilizer subgroup $S \subset \PGLn$ 
is reductive, and hence so is the preimage $S^*$ of $S$ in $\GL_n$. 
Thus the natural $n$-dimensional linear representation 
$\phi \colon S^* \hookrightarrow \GL_n$ decomposes
as a direct sum of irreducibles. Suppose $\phi$ has
irreducible components $\phi_1, \dots, \phi_r$
of dimensions $d_1, \dots, d_r$, occurring with
multiplicities $e_1, \dots, e_r$, respectively.  
Note that $d_1 e_1 + \dots + d_r e_r = n$. Our main result is the following
proposition. Corollary~\ref{cor.semisimple} is an immediate consequence  
of parts (a), (b) and (c).

 \begin{prop} \label{prop.semisimple} 
\begin{enumerate}
 \item[(a)]$A^G$ is a form of $\Mat_{e_1}(k) \times \dots \times
   \Mat_{e_r}(k)$.  In particular, $A^G$ is semisimple.
 \item[(b)]The maximal degree of an element of $A^G$ over $K$ is $e_1 +
   \dots + e_r$.
 \item[(c)]The maximal degree of an element of $\Z(A^G)$ over $K$ is
  $r$.
 \item[(d)] If $A^G$ is a simple algebra, then the normalizer $N(S^*)$
   of $S^*$ in $\GL_n$ transitively permutes  
   {\upshape(}the isomorphism types of\/{\upshape)} the irreducible
   representations $\phi_1, 
   \dots, \phi_r$ of $S^*$ occurring in $\phi$.  In particular, 
   in this case, $d_1 = \dots = d_r$ and $e_1 = \dots = e_r$.
 \end{enumerate}
\end{prop}

\begin{proof}
  (a) By a corollary of Schur's Lemma (see, e.g.,
  \cite[Proposition~1.8]{farb-dennis}),
  \[ \Mn^S=(\Mn)^{S^*}
  =\operatorname{End}_{S^*}(k^n) \simeq
  \Mat_{e_1}(k) \times \dots \times \Mat_{e_r}(k) \, . \]
Part (a) now follows from Theorem~\ref{thm2}.

  \smallskip
(b) The degree of $x\in A^G$ over $K$ is $\dim_K(K[x])$.
In view of part (a), the maximal degree of an element of $A^G$ 
over $K$ is the same as the maximal degree of an element 
of $\Mat_{e_1}(k) \times \dots \times \Mat_{e_r}(k)$ over $k$.
The latter is easily seen to be $e_1 + \dots + e_r$.

\smallskip
(c) Once again, in view of part (a), we only need to show that the
maximal degree of a central element of
$\Mat_{e_1}(k) \times \dots \times \Mat_{e_r}(k)$ over $k$
is $r$. The center of $\Mat_{e_1}(k) \times \dots \times \Mat_{e_r}(k)$
is isomorphic, as a $k$-algebra, to $k \times \dots \times k$ ($r$ times),
and part (c) follows.

\smallskip
 (d) Let $V_i$ be the unique $S^*$-subrepresentation 
of $k^n$ isomorphic to $\phi_i^{e_i}$. The normalizer $N(S^*)$ 
clearly permutes the subspaces $V_i$ (and hence, the representations
$\phi_i$); we want to show that this permutation
action is transitive. Assume the contrary: say,
$\{ \phi_1, \dots, \phi_s \}$ is an $N(S^*)$-orbit 
for some $s < r$. Then $W_1 = V_1 \oplus \dots \oplus V_s$ and
$W_2 = V_{s+1} \oplus \dots \oplus V_r$ are complementary 
$N(S^*)$-invariant subspaces of $k^n$ and
\[ \Mn^S=\operatorname{End}_{S^*}(k^n) \simeq
\operatorname{End}_{S^*}(W_1) \oplus 
\operatorname{End}_{S^*}(W_2) \, , \]
where $\simeq$ denotes an $N(S)$-equivariant isomorphism of
$k$-algebras.
Now let $X$ be an associated $G\times\PGLn$-variety for the action 
of $G$ on $A$, $Y = X/G$, and $Z$ be the union of the components 
of maximal dimension
  of $Y^S$, as in the proof of Proposition~\ref{prop.thm2}.  Then,
  \begin{align*}
  A^G &\stackrel{\text{\tiny \eqref{e.X->Y}}}{=}
  \RMaps_{\PGLn}(Y, \Mn)  
  \stackrel{\text{\tiny \eqref{e.Y->Zrevived},\eqref{e.Z}}}{=} 
  \RMaps_{N(S)}(Z, \Mn^S)  \\
  &\stackrel{\hphantom{\text{\tiny \eqref{e.X->Y}}}}{=}  
  \RMaps_{N(S)}(Z, \operatorname{End}_{S^*}(W_1)) \oplus 
\RMaps_{N(S)}(Z, \operatorname{End}_{S^*}(W_2)) \, . 
  \end{align*}
  This shows that $A^G$ is not simple.
\end{proof}

\section{Proof of Theorem~\ref{thm.gal-cohomology}}
\label{sect.thm.gal-cohomology}

We begin with some preliminaries on Galois cohomology. 
Let $F/k$ be a finitely generated field extension and 
$\Gamma/k$ be a linear algebraic group. Let $\Gamma$-$\Var(F)$ be the
set of isomorphism classes of 
generically free $\Gamma$-varieties $X$ with $k(X)^\Gamma = F$. 
Note that since $k(X)^\Gamma$ is a field, $X$ is necessarily primitive,
i.e., $\Gamma$ transitively permutes the irreducible components 
of $X$; see Section~\ref{sect.fiber-spaces}.
Recall that $H^1(F, \Gamma)$ is in a natural bijective correspondence
with $\Gamma$-$\Var(F)$; see, e.g.,~\cite[(1.3)]{popov}. 
We can thus identify $H^1(F, \Gamma)$ 
with $\Gamma$-$\Var(F)$. If $\alpha \colon \Gamma \to \Gamma'$ is 
a homomorphism of algebraic groups then the induced map 
$\alpha_* \colon H^1(F, \Gamma) \to H^1(F, \Gamma')$ is given by
$Z \mapsto \Gamma' *_\Gamma Z$; this follows from
\cite[Theorem 1.3.3(b)]{popov} (see also~\cite[Proposition 28.16]{boi}).

If $\Gamma$ is the automorphism group of some finite-dimensional
$k$-algebra $R$ then, as we remarked before the statement of
Theorem~\ref{thm.gal-cohomology}, elements of 
$H^1(F, \Gamma)$ are also in a natural bijective correspondence 
with the set of isomorphism classes of $F$-forms of $R$.
These two interpretations of $H^1(F, \Gamma)$ can be related as follows:
the algebra corresponding to the $\Gamma$-variety $W$ is 
$B = \RMaps_\Gamma(W, R)$; cf.~\cite[Proposition 8.6]{reichstein}.  
Here, as usual, $\RMaps_\Gamma$ stands for $\Gamma$-equivariant 
rational maps.

We are now ready to proceed with the proof 
of Theorem~\ref{thm.gal-cohomology}. We fix $n \ge 1$ and a subgroup 
$S \subset \PGLn$. For notational convenience, set 
$H = \Aut(\Mn^S)$. Consider the following diagram:
\[ \xymatrix@R-18pt{ H^1(F, N/S) \ar@{->}[r]^-{\alpha_*} \ar@{=}[d] & 
              H^1(F, H)  \ar@{=}[d] &        \cr
       \text{$N/S$-$\Var(F)$} \ar@{->}[r] & 
\text{$H$-$\Var(F)$} \ar@{->}[r]^-{\simeq} &  
\{\text{$F$-forms of $\Mn^S$}\} \cr
   Z \ar@{|->}[r] & H*_{N/S} Z^{\rule{0pt}{1.5ex}} &       \cr 
    & W \ar@{|->}[r] &  \RMaps_H(W, \Mn^S) }
\]
\smallskip

\noindent
This diagram shows that the image of the map 
$\alpha_* \colon H^1(F, N/S) \to H^1(F, H)$
consists of all classes in $H^1(F, H)$ 
corresponding to algebras $B/F$ of the form 
\[ B = \RMaps_H(H *_{N/S} Z, \Mn^S) \, , \]
where $F = k(H*_{N/S} Z)^H$ 
and $Z$ is a primitive $N/S$-variety. But
\[ \RMaps_H(H *_{N/S} Z, \Mn^S)
    \simeq \RMaps_{N/S} (Z, \Mn^S)\]
and $k(H*_{N/S} Z)^H \simeq k(Z)^{N/S}$;
see Lemma~\ref{lem.Y->Z}. 
Consequently, $[B]$ lies in the image of $\alpha_*$ if 
and only if  there exists a generically free primitive
$N/S$-variety $Z$ such that $B \simeq \RMaps_{N/S}(Z, \Mn^S)$ and
$F \simeq k(Z)^{N/S}$, i.e., if and only if condition (c) of the
following lemma is satisfied.

We will say that an $F$-algebra $B$ is a {\em fixed algebra with
  associated stabilizer} $S \subset \PGLn$ if there exists a tame
action of an algebraic group $G$ on a central simple algebra $A/K$ of
degree $n$ with associated stabilizer $S$ such that $B \simeq A^G$ and
$F \simeq K^G$. In other words, $B/F$ is a fixed algebra with associated
stabilizer $S$ if it satisfies condition (a) of
Theorem~\ref{thm.gal-cohomology}.
Hence Theorem~\ref{thm.gal-cohomology} follows from the following result:

\begin{lem} \label{lem1.gal-cohomology} Let $B$ be an $F$-algebra,
  $S$ a subgroup of $\PGLn$, and $N=N_\PGLn(S)$.  The following are
  equivalent:

\begin{enumerate}
\item[(a)]$B/F$ is a fixed algebra with associated stabilizer $S$.
\item[(b)]There exists an irreducible $\PGLn$-variety $Y$ with
  stabilizer $S$ in general position such that $B \simeq
  \RMaps_{\PGLn}(Y, \Mn)$ and $F \simeq k(Y)^{\PGLn}$.
\item[(c)]There exists a generically free primitive $N/S$-variety $Z$
  such that $B \simeq \RMaps_{N/S}(Z, \Mn^S)$ and $F \simeq k(Z)^{N/S}$.
\end{enumerate}
\end{lem}

The isomorphisms in the various parts of the lemma are compatible: the
second is always a restriction of the first.

\begin{proof} (a) $\Longrightarrow$ (b):
By definition the algebras $B/F$ satisfying condition (a) 
of Theorem~\ref{thm.gal-cohomology} are those of the form 
$B = \RMaps_{\PGLn}(X, \Mn)^G$,
where $F = k(X)^{G \times \PGLn}$. Here $X$ is an irreducible $G
\times \PGLn$-variety, 
such that the $\PGLn$-action on $X$ is generically free,
and the projection of $\Stab_{G \times \PGLn}(x)$ to $\PGLn$ 
is conjugate to $S$ for $x \in X$ in general position. 
To write $B$ and $F$ as in (b), take $Y = X/G$; see~\eqref{e.X->Y}.

(b) $\Longrightarrow$ (c): Take $Z$ to be the union 
of the components of $Y^S$ of maximal dimension, as in Step~2 of the
proof of Proposition~\ref{prop.thm2}.

(c) $\Longrightarrow$ (b): Take $Y$ to be the homogeneous fiber space
$\PGLn *_N Z$, and use Lemma~\ref{lem.Y->Z} and~\eqref{e.Z}.
Note that $Y$ is irreducible since $k(Y)^\PGLn$ is a field and $\PGLn$
is connected.  A simple computation shows that $S$ is the stabilizer
in general position for the $\PGLn$-action on $Y$.

(b) $\Longrightarrow$ (a): Use Lemma~\ref{lem.wild}, with
$H = \PGLn$ to reconstruct $X$ from $Y$. Note that here $G$ 
can be any connected linear algebraic group containing $\PGLn$.

This completes the proof of Lemma~\ref{lem1.gal-cohomology} 
and thus of Theorem~\ref{thm.gal-cohomology}.
\end{proof}

\section{An example}

We will now illustrate Theorem~\ref{thm.gal-cohomology} with an example.
Fix an integer $n \ge 1$ and write $n = d_1 e_1 + \dots + d_r e_r$
for some positive integers $d_1, e_1, \dots, d_r, e_r$. 
A partition of $n$ of this form naturally arises 
when an $n$-dimensional representation is decomposed 
as a direct sum of irreducibles of dimensions
$d_1, \dots, d_r$ occurring with multiplicities 
$e_1, \dots, e_r$, respectively; in particular, the integers $d_i$ 
and $e_i$ came up in this way at the beginning of 
Section~\ref{sect.very-tame}.
For this reason $\tau = [(d_1, e_1), \dots, (d_r, e_r)]$ 
is often referred to as a ``representation type''; cf.~\cite{lbp}.

Choose vector spaces $V_1, W_1, \dots, V_r, W_r$, so that
$\dim V_i = d_i$ and $\dim W_i = e_i$ and an isomorphism
\[ k^n \simeq (V_1 \otimes W_1 )\oplus \dots \oplus (V_r \otimes W_r) \, . \]
Let $S_{\tau}^* \simeq \GL_{d_1} \times \dots \times \GL_{d_r}$ 
be the subgroup
\[ S_{\tau}^* = (\GL(V_1) \otimes I_{W_1}) \times \dots \times
(\GL(V_r) \otimes I_{W_r}) \subset \GL_n \, , \]
where $I_{W}$ denotes the identity linear transformation $W \to W$.
Let $N_{\tau}^*$ be the normalizer of $S_{\tau}^*$ in $\GL_n$.
Denote the image of  $S_{\tau}^*$ under the natural projection 
$\GL_n \to \PGLn$ by $S_{\tau}$ and the normalizer of $S_{\tau}$ in
$\PGLn$ by $N_{\tau}$. The purpose of this section is 
to spell out what Theorem~\ref{thm.gal-cohomology} tells us for 
$S = S_{\tau} \subset \PGLn$. 
Note that
\[ \Mn^{S_{\tau}} = (I_{V_1} \otimes \End(W_1)) \times \dots \times 
(I_{V_r} \otimes \End(W_r)) \simeq \Mat_{e_1} \times \dots \times 
\Mat_{e_r} \, . \] 
Our goal is to determine which forms of $\Mn^{S_{\tau}}$
occur as fixed algebras with associated stabilizer $S_{\tau}$.

Note that reordering the pairs $(d_i, e_i)$ 
will have the effect of replacing
$S_{\tau}$ by a conjugate subgroup of $\PGLn$. Since the group
$S$ in Theorems~\ref{thm2} and \ref{thm.gal-cohomology} is only
defined up to conjugacy, we are free to order them any
way we want. In particular, if there are $t$ distinct numbers
$f_1, \dots, f_t$ among $e_1, \dots, e_r$, occurring with 
multiplicities $r_1, \dots, r_t$, respectively, we will, for
notational convenience, renumber the pairs $(d_i, e_i)$
so that $\tau$ assumes the form               
\[ \tau = [(d_{11}, f_1), \dots, (d_{1r_1}, f_1), \dots, (d_{t1}, f_t),
\dots, (d_{t r_t}, f_t)] \, . \]               
Suppose that for each fixed $i$, exactly $l_i$ numbers occur among
$d_{i1}, \dots, d_{ir_i}$, with multiplicities $r_{i1}, r_{i2}, \dots,
r_{il_i}$. We order the pairs $(d_{i1}, f_i), \dots, (d_{ir_i}, f_i)$
so that the first $r_{i1}$ of the $d_{ij}$s are the same, the next
$r_{i2}$ are the same, and so on.  Clearly $ r_i = r_{i1} + \dots +
r_{i l_i}$,
\[ \Mn^{S_{\tau}} \simeq \Mat_{e_1} \times \dots \times \Mat_{e_r}  
 \simeq \Mat_{f_1}^{r_1} \times \dots \times \Mat_{f_t}^{r_t} \]
and
\[ \Aut(\Mn^{S_{\tau}}) \simeq \prod_{i = 1}^t \Aut(\Mat_{f_i}^{r_i})  \, . \]
Note that $\Aut(\Mat_f^r)$ is the semidirect product 
$\PGL_f^r \sdp \Sym_r$, where $\PGL_f^r$ is the group of
inner automorphisms of $\Mat_f^r$, and $\Sym_r$ acts by permuting 
the $r$ factors of $\Mat_f$.

\begin{lem} \label{lem.alpha}
The natural homomorphism $\alpha\colon N_{\tau}/S_{\tau} \to 
\Aut(\Mn^{S_{\tau}})$ corresponds to the inclusion map
\[ \prod_{i= 1}^t \biggl( \prod_{j=1}^{l_i} \Aut(\Mat_{f_i}^{r_{ij}}) \biggr)
 \hookrightarrow
\prod_{i=1}^t \Aut(\Mat_{f_i}^{r_{i}}) \, . \]
{\upshape(}Recall that $r_i = r_{i1} + \dots + r_{i l_i}$ for each
$i$.{\upshape)} 
\end{lem}

Note that this lemma says that $\alpha$
embeds $N_{\tau}/S_{\tau}$ in $\Aut(\Mn^{S_{\tau}})$ as 
a subgroup of finite index. In particular, the image of $\alpha$
contains the connected component
\[ \PGL_{e_1} \times \cdots \times \PGL_{e_r}
      = \PGL_{f_1}^{r_1} \times \dots \times \PGL_{f_t}^{r_t} \]
of $\Aut(\Mn^{S_{\tau}})$.

\begin{proof} First of all, note that we may replace 
$N_{\tau}/S_{\tau}$ by $N_{\tau}^*/S_{\tau}^*$; these two groups
are isomorphic and act the same way on $\Mn^{S_{\tau}}$.

We claim that the homomorphism $\alpha \colon N_{\tau}^*/S_{\tau}^* \to 
\Aut(\Mn^{S_{\tau}})$ is injective. Indeed, the 
kernel of the action of $N_{\tau}^*$ on $\Mn^{S_{\tau}}$ 
is the centralizer of 
\[ \Mn^{S_{\tau}} = (I_{V_1} \otimes \End(W_1)) \times \dots \times 
(I_{V_r} \otimes \End(W_r)) \]
in $\GL_n$. Using Schur's Lemma, as in the proof of 
Proposition~\ref{prop.semisimple}(a), we see that this centralizer is  
\[ S_{\tau}^* = 
(\GL(V_1) \otimes I_{W_1}) \times \dots \times 
(\GL(V_r) \otimes I_{W_r} )\, , \]
and the claim follows.  Thus $\alpha$ is injective, and we 
only need to determine which automorphisms of
$\Aut(\Mn^{S_{\tau}}) \simeq \Mat_{e_1} \times \dots \times \Mat_{e_r}$ 
can be realized as conjugation by some $g \in N_{\tau}^*$. 
Taking 
$g = (I_{V_1} \otimes a_1, \dots , I_{V_r} \otimes a_r)$ 
and letting $a_i$ range over $\GL(W_i)$,
we see that every element of 
\[ \PGL_{e_1} \times \cdots \times \PGL_{e_r}
      = \PGL_{f_1}^{r_1} \times \dots \times \PGL_{f_t}^{r_t} \]
lies in the image of $\alpha$. 

Let us now view $k^n$ as an $n$-dimensional
representation of the group $S_{\tau}^*\simeq\GL_{d_1} \times \dots
\times \GL_{d_r}$, as we did at the beginning of this section.
Any $g \in N_{\tau}^*$ permutes the isomorphism types of the
irreducible representations
$V_1, \dots, V_r$ of $S_{\tau}^*$; it also permutes the isotypical 
components $V_1 \otimes W_1, \dots, V_r \otimes W_r$ of this representation
(as subspaces of $k^n$).
If $g$ maps the isotypical component
$V_i \otimes W_i$ to the isotypical component $V_j \otimes W_j$
then their dimensions $d_ie_i$ and $d_j e_j$ have to be the same.
The dimensions $d_i$ and $d_j$ of the underlying irreducible 
representations $V_i$ and $V_j$ also have to be the same.
In other words, this is only possible if $(d_i, e_i) = (d_j, e_j)$. 

Conversely, if $(d_i, e_i) = (d_j, e_j)$ then a suitably chosen permutation 
matrix $g \in N_{\tau}^*$, will interchange
$V_i \otimes W_i$ and $V_j \otimes W_j$ and preserve every 
other isotypical component $V_m \otimes W_m$ for $m \ne i, j$.  
These permutation matrices generate the finite subgroup 
$\prod_{i= 1}^t \prod_{j=1}^{l_i} \Sym_{r_{ij}}$ of 
$\Aut(\Mn^{S_{\tau}})$.

Now let $g$ be any element of $N_{\tau}^*$.  After composing $g$ 
with a product of permutation matrices as above,
we arrive at a $g_0 \in N_{\tau}^*$ which preserves 
every isotypical component $V_i \otimes W_i$. Conjugating 
by such a $g_0$ will preserve every factor
$I_{V_i} \otimes \End(W_i) \simeq \Mat_{e_i}$ 
of $\Aut(\Mn^{S_{\tau}})$ and thus will lie in
$\PGL_{f_1}^{r_1} \times \dots \times \PGL_{f_t}^{r_t}$. 
We conclude that the image of $\alpha$ is generated by
\[ \text{$\PGL_{f_1}^{r_1} \times \dots \times \PGL_{f_t}^{r_t}$ 
and
$\prod_{i= 1}^t \prod_{j=1}^{l_i} \Sym_{r_{ij}}$}  \]
and the lemma follows.
\end{proof}

\begin{prop} \label{prop.tau}
$B/F$ is a fixed algebra with associated stabilizer $S_{\tau}$ 
if and only if $B$ is $F$-isomorphic to the direct product
\[ (B_{11} \times \dots \times B_{1l_1}) \times \dots \times (B_{t1} \times 
\dots \times B_{t l_t}) \]
where $B_{ij}/F$ is a form of $\M_{f_i}^{r_{ij}}/k$.
\end{prop}

\begin{proof} Recall that
$H^1(F, \Aut(R)) = \{ \text{$F$-forms of $R$} \}$. 
Now suppose $R_1$, \ldots, $R_m$ are $k$-algebras,
$R = R_1 \times \dots \times R_m$, $G_i = \Aut(R_i)$ and
$G = \Aut(R_1) \times \dots \times \Aut(R_m)$. Then the natural 
inclusion
$G_1 \times \dots \times G_r \hookrightarrow \Aut(R)$
induces a morphism
from
\[ H^1(F, G_m \times \dots \times G_m) \simeq 
 H^1(F, G_1) \times \dots \times H^1(F, G_m)\]
to $H^1(F, \Aut(R))$
which maps $(A_1, \dots, A_m)$ to $A = A_1 \times \dots \times A_m$.
Here $A_i$ is an $F$-form of $R_i$ (and hence, $R$ is an $F$-form of $A$).

In the language of $G$-varieties, the first isomorphism says 
that $X$ is a generically free primitive $G$-variety if and only
if $X$ is birationally isomorphic (as a $G$-variety) to the fiber 
product $X_1 \times_Y \dots \times_Y X_m$, where $X_i$ is 
a primitive generically free $G_i$-variety such that $X_i/G_i \simeq Y$
is a model for $F$ ($Y$ is the same for each $i$).
Note that $X_1$ can be recovered from $X$
as the rational quotient $X/(G_2 \times \dots \times G_m)$; similarly for 
$X_2, \dots, X_m$.  The second map is given by the natural isomorphism
\[ \RMaps_G(X, R) \simeq 
\RMaps_{G_1}(X_1, R_1) \times \dots \times \RMaps_{G_m}(X_m, R_m)  \] 
of $F$-algebras.

We now apply this to the homomorphism $\alpha \colon N_{\tau}/S_{\tau} \to
\Aut(\Mn^{S_{\tau}})$ described in Lemma~\ref{lem.alpha}, and 
the proposition follows from Theorem~\ref{thm.gal-cohomology}.
\end{proof}

We conclude with two corollaries of Proposition~\ref{prop.tau}.

\begin{cor} \label{cor2.tau}
  Let $\tau = [(d_1, e_1), \dots, (d_r, e_r)]$.  Then the following
  are equivalent.
  \begin{enumerate}
  \item[(a)] $d_i = d_j$ whenever $e_i = e_j$.
  \item[(b)] Every form $B/F$ of $\Mn^{S_{\tau}}$ appears as a fixed
    algebra with associated stabilizer $S_{\tau}$.
  \end{enumerate}
\end{cor}

\begin{proof} (a) is equivalent to $l_i = 1$ for every $i = 1, \dots, t$. 
The implication (a) $\Rightarrow$ (b) is now an immediate consequence of
the proposition.  To show that (b) $\Longrightarrow$ (a),
consider an $F$-algebra $B = D_1 \times \dots \times D_t$, where
$D_i$ is a division algebra of degree $f_i$ such that
$[\Z(D_i): F] = r_i$. Here $\Z(D_i)$ denotes the center of $D_i$.
For example, we can take $F = k(a_1, b_1, \dots, a_t, b_t, c)$, 
where $a_1, b_2, \dots, a_t, b_t, c$ are independent variables,
and $D_i$ is the symbol algebra $(a_i, b_i)_{f_i}$ over 
$\Z(D_i) = F(\root{r_i} \; \; \of c)$.
The algebra $B = D_1 \times \dots \times D_t$ is then a form of
\[ \Mn^{S_{\tau}} = \Mat_{f_1}^{r_1} \times \dots \times \M_{f_t}^{r_t} \, . \]
On the other hand, if (b) holds, Proposition~\ref{prop.tau} tells us that
$B$ can be written as a direct product of $l_1 + \dots + l_t$
$F$-algebras. But $B = D_1 \times \dots \times D_t$ clearly
cannot be written 
as a non-trivial direct product of more than $t$ $F$-algebras. 
Thus $l_1 + \dots + l_t \le t$, i.e., $l_1 = \dots = l_t = 1$, 
and (a) follows.
\end{proof}

\begin{cor} \label{cor3}
  Every finite-dimensional semisimple algebra $B/F$ {\upshape(}where
  $F$ is a finitely generated field extension of $k${\upshape)}
  appears as $A^G/K^G$ for some very tame geometric action of a linear
  algebraic group $G$ on a central simple algebra $A$ with center $K$.
\end{cor} 

Note that here we make no assumption on the degree of $A$ or on
the associated stabilizer. However, the proof will show that 
the associated stabilizer can always be taken to be the reductive
group $S_{\tau}$, where
$\tau = [(1, e_1), \dots, (1, e_r)]$ for suitable positive integers 
$e_1, \dots, e_r$.

\begin{proof}
  Every semisimple algebra is a form of $\cA = \M_{e_1} \times \dots
  \times \Mat_{e_r}$ for some positive integers $e_1, \dots, e_r$. Now
  set $n = e_1 + \dots + e_r$, $d_1 = \dots = d_r = 1$ and $\tau =
  [(1, e_1), \dots, (1, e_r)]$.  Then $\Mn^{S_{\tau}} \simeq \cA$ and
  the conditions of part (a) of Corollary~\ref{cor2.tau} are satisfied
  for this choice of $\tau$. Hence, every form of $\cA$ appears as a
  fixed algebra (with associated stabilizer $S_{\tau}$).
\end{proof}

\end{document}